\documentclass[12pt]{article}
\usepackage[T1]{fontenc}
\usepackage{lmodern}
\usepackage{amsmath,amssymb,amsthm,mathtools}
\usepackage[margin=1in]{geometry}
\usepackage{microtype}
\usepackage{enumitem}
\usepackage[dvipsnames]{xcolor}
\usepackage{tikz}
\usetikzlibrary{calc}
\usepackage[colorlinks=true,linkcolor=MidnightBlue,citecolor=MidnightBlue,urlcolor=MidnightBlue]{hyperref}
\newtheorem{theorem}{Theorem}[section]
\newtheorem{prop}[theorem]{Proposition}
\newtheorem{lemma}[theorem]{Lemma}
\newtheorem{cor}[theorem]{Corollary}
\newtheorem{claim}[theorem]{Claim}
\newtheorem{obs}[theorem]{Observation}
\newtheorem{question}[theorem]{Question}
\theoremstyle{definition}
\newtheorem{defn}[theorem]{Definition}
\newtheorem{remark}[theorem]{Remark}
\renewcommand{\leq}{\leqslant}
\renewcommand{\le}{\leqslant}
\renewcommand{\geq}{\geqslant}
\renewcommand{\ge}{\geqslant}
\newcommand{\defeq}{\vcentcolon=}

\newcommand{\PP}{\mathbb P}
\newcommand{\ZZ}{\mathbb Z}
\DeclareMathOperator{\Cay}{Cay}
\DeclareMathOperator{\mis}{mis}
\newcommand{\triangledef}{\operatorname{def}_{\triangle}}
\setlist[enumerate]{leftmargin=*,itemsep=2pt,topsep=4pt}
\title{On the number of tiles in $[n]^d$}
\author{Anubhab Ghosal\thanks{Mathematical Institute, University of Oxford, Oxford OX2 6GG, UK. Email: \texttt{ghosal@maths.ox.ac.uk}. Supported by the Clarendon Fund and Oxford Ryniker Lloyd Graduate Scholarship.} \and Jeck Lim\thanks{Institute for Applied Mathematics, University of Bonn, 53115 Bonn, Germany. Email: \texttt{jlim@uni-bonn.de}. Supported by the Deutsche Forschungsgemeinschaft (DFG, German Research Foundation) under Germany's Excellence Strategy -- EXC-2047/2 -- 390685813. Part of this work was done while the author was at the University of Oxford supported by ERC Advanced Grant 883810.} \and Elad Tzalik\thanks{Weizmann Institute of Science, Rehovot, Israel.
Email: \texttt{Elad.tzalik@weizmann.ac.il}.
Supported by the Adams Fellowship Program of the Israel Academy of Sciences and Humanities.}}
\date{}
\begin{document}
\maketitle

\begin{abstract}
 Let $\mathcal{T}(n,d)$ denote the set of subsets of $[n]^d$ that tile $\mathbb Z^d$ by translations, and write $T(n,d)=|\mathcal{T}(n,d)|$. We show that almost every tile in $\mathcal{T}(n,1)$ tiles $\ZZ$ by an arithmetic progression, thus determining $T(n,1)$ asymptotically. In higher dimensions, we prove that $T(n,d)=O_d\!\left(n^{C_d}3^{n^d/3}\right)$, confirming a conjecture of Benjamini, Kozma and Tzalik. In the other direction, we prove that $T(n,d)\geq (1-o(1))n^{\binom d2}3^{n^d/3}$ when $3\mid n$. This shows one cannot take $C_d=d-1$ for $d\geq 3$, answering a question due to the same authors negatively. Finally, we also obtain nearly tight bounds in dimension two. Our upper bounds rely on a new observation realising tiles as maximal independent sets in truncated Cayley graphs.
\end{abstract}

{
\par
}


\section{Introduction}

We study the number and typical structure of finite sets that tile
$\ZZ^d$ by translations.

\begin{defn}
For sets $A,Z\subseteq\ZZ^d$, we write $A\oplus Z$ for the sumset $A+Z$,
together with the assertion that every $v\in A+Z$ has a unique representation
$v=a+z$ with $a\in A$ and $z\in Z$. We say that $A$ \emph{tiles $\ZZ^d$
by $Z$} if $\ZZ^d=A\oplus Z$. We call $A$ a \emph{tile}
if such a set $Z$ exists.
\end{defn}

The structure of finite tiles has attracted sustained interest from
mathematicians for decades, connecting number theory, combinatorics,
and harmonic analysis. Even in dimension one, characterising all finite tiles
has been a longstanding problem. In a foundational work in 1977,
Newman~\cite{Newman1977} characterised tiles of prime-power cardinality.
Subsequently, in 1999, Coven and Meyerowitz~\cite{CovenMeyerowitz1999} gave arithmetic
conditions that are sufficient for a finite set to tile $\ZZ$, and conjectured that these conditions are also necessary. Coven and Meyerowitz proved their conjecture in the case that the tile's cardinality has at most two distinct prime factors. Extending this to even three prime factors proved to be surprisingly difficult, and only recently,
{\L}aba and
Londner~\cite{LabaLondner2023} were able to resolve particular three-prime cases.
In a very recent preprint, Prajapati~\cite{Prajapati2026} claims a full resolution
of the conjecture, accompanied
by a Lean formalisation.

In higher dimensions, tiles exhibit increasingly complex behaviour. Consider, for example, the periodic tiling conjecture. A tile in $\ZZ^d$ tiles periodically if it tiles $\ZZ^d$ by a finite union of cosets of a finite-index subgroup of $\ZZ^d$. The celebrated periodic tiling conjecture~\cite{Stein1974,GrunbaumShephard1987,LagariasWang1996} asserts that all finite tiles in $\ZZ^d$ tile periodically. As first observed by Newman~\cite{Newman1977}, a simple pigeonhole argument shows that the periodic tiling conjecture holds when $d=1$. For $d=2$, the conjecture was established by Bhattacharya~\cite{Bhattacharya2020} by combining techniques from ergodic theory and number theory. In a substantial breakthrough, Greenfeld and Tao~\cite{GreenfeldTao2024}
disproved the periodic tiling conjecture in sufficiently high dimension by constructing a \emph{Sudoku puzzle} that can be encoded as a tiling equation and has necessarily aperiodic solutions.
In a preprint that appeared very recently (with an accompanying Lean verification), OpenAI~\cite{OpenAI2026} claims to extend Greenfeld and Tao's work to construct a counterexample to the periodic tiling conjecture in dimension $d=3$. For an extended discussion on \emph{wild} tilings in high dimensions, see Greenfeld's survey~\cite{Greenfeld2026}.

These difficulties motivate the question: what does a
typical tile in $\ZZ^d$ look like? In spite of the complexity of tilings, it is possible that \emph{almost all} tiles afford simpler descriptions.
Benjamini, Kozma, and Tzalik~\cite{BenjaminiKozmaTzalik2026}
pursued this viewpoint by counting tiles contained in a finite box.
\begin{defn}
We define $T(n,d)$ to be the number of sets
$A\subseteq[n]^d$ that tile $\ZZ^d$.
\end{defn}
\noindent Benjamini, Kozma and Tzalik~\cite{BenjaminiKozmaTzalik2026}
proved that, for every fixed $d\ge1$,
\begin{equation}
\label{BKT bound}
    T(n,d)=3^{n^d/3+ o(n^d)}.
\end{equation}
They further conjectured that the upper
bound could be sharpened to a polynomial in $n$ times $3^{n^d/3}$.
Our first result confirms this conjecture.

\begin{theorem}\label{thm:general-upper}
For every fixed $d\ge1$, there is a constant $C_d$ such that
$T(n,d)=O_d(n^{C_d}3^{n^d/3})$.
\end{theorem}

Additionally, in dimension one, we obtain a structural description of a typical tile
and precise asymptotics for the count.

\begin{defn}
We call a tile $A\subseteq\ZZ$ an \emph{AP tile} if
$\ZZ=A\oplus k\ZZ$ for some integer $k\ge1$, and a \emph{non-AP tile}
if no such $k$ exists. Let $T_{\mathrm{AP}}(n,1)$ and
$T_{\mathrm{non\text{-}AP}}(n,1)$ denote the numbers of AP and non-AP tiles
contained in $[n]$, respectively.
\end{defn}

\noindent AP tiles are precisely the sets which tile $\ZZ$ by an arithmetic progression, hence the name. In dimension one, we show that almost every tile contained in
$[n]$ is an AP tile, which also implies an asymptotic formula for the number of tiles contained in $[n]$.

\begin{theorem}\label{thm:one-dimensional}
We have $T_{\mathrm{non\text{-}AP}}(n,1)=o(3^{n/3})$. Moreover, if $n=3m+s$
with $s\in\{0,1,2\}$, then
\[
    T(n,1)=(c_s+o(1))3^m,
    \qquad
    (c_0,c_1,c_2)=\left(\frac{217}{17},\frac{312}{17},\frac{450}{17}\right).
\]
\end{theorem}

In dimension two, a more involved argument gives an upper
bound that matches the lower bound (in Theorem~\ref{thm:general-lower}) up to polylogarithmic (in $n$)
factors when $3\mid n$.

\begin{theorem}\label{thm:two-dimensional}
We have $T(n,2)=\widetilde O(n3^{n^2/3})$.
\end{theorem}
\noindent Here, and in the rest of the paper, $\widetilde O(g(n))$ is shorthand for $O(\log^{O(1)}(n)g(n))$.

Benjamini, Kozma and Tzalik~\cite{BenjaminiKozmaTzalik2026} proved the upper bound in Equation~\eqref{BKT bound} using an entropy-based approach. Our upper bounds build on a new observation connecting tilings and maximal independent sets in truncated Cayley graphs. For a prescribed set of
translations, each tile is a maximal independent set in an associated truncated Cayley graph; see Lemma~\ref{lem:cover-mis} for a precise formulation. Moon and Moser~\cite{MoonMoser1965} proved that an $N$-vertex
graph has at most $3^{N/3}$ maximal independent sets. Applying the Moon-Moser bound to these Cayley graphs, one obtains a proof of Theorem~\ref{thm:general-upper} almost immediately.

In dimensions one and two, we use the stability
version of the Moon-Moser theorem, due to Kahn and Park~\cite{KahnPark2020StabilityMIS},
to show that graphs contributing substantially to the count must
contain a large induced subgraph consisting of disjoint triangles.
We then classify the translation configurations that can produce
this structure. This is relatively straightforward in dimension one, but owing to richer geometry with increase in  dimension, the analysis becomes more technical in dimension two.

Benjamini, Kozma and Tzalik~\cite{BenjaminiKozmaTzalik2023} also proved a lower bound of $n^{d-1}\cdot 3^{n^d/3}$ for $T(n,d)$ when $3 \mid n$. Motivated by this, they asked if one can take $C_d=d-1$ in the upper bound in Theorem~\ref{thm:general-upper}. We improve their lower bound, answering the question negatively when $d\geq 3$. 

\begin{theorem}\label{thm:general-lower}
For every fixed $d\ge1$, if $3\mid n$, then
\[
    T(n,d)\ge (1-o(1))n^{\binom d2}3^{n^d/3}.
\]
\end{theorem}

\noindent Our strategy for proving the lower bound is to construct $n^{\binom d2}$ lattices, so that for each lattice $L$, there is a family of $3^{n^d/3}$ tiles in $[n]^d$ which tile $\mathbb{Z}^d$ by $L$, and these families are mostly disjoint for distinct lattices.

\paragraph{Organisation.}The rest of the paper is organised as follows. We prove Theorems~\ref{thm:general-upper}, \ref{thm:one-dimensional}, \ref{thm:two-dimensional} and \ref{thm:general-lower} in Sections~\ref{sec:upper-all}, \ref{sec:upper-one}, \ref{sec:upper-two} and \ref{sec:lower} respectively. Finally, in Section~\ref{sec: conc}, we highlight some open problems and further directions for research.

\paragraph{Use of AI.}
All of the main ideas in this paper were developed by the authors without any AI assistance. We used ChatGPT 6 Pro to assist substantially with writing and exposition, and ChatGPT 5.4 Pro to
calculate the constants $c_s$ in Theorem~\ref{thm:one-dimensional}.

\paragraph{Acknowledgements.}
A.G. and E.T. thank IAS/PCMI for organising the 2025 Graduate
Summer School on Extremal and Probabilistic Combinatorics, where this collaboration was initiated. In addition, A.G. thanks IAS/PCMI for funding his visit to the summer school. E.T. thanks Itai Benjamini for his encouragement over the course of this project. A.G. thanks Emma Hogan for coining the term \emph{crossant}.

\section{Upper bound in all dimensions}\label{sec:upper-all}
The main observation behind our upper bound is that, for a prescribed translation
set $Z$, the number of tiles is bounded by the number of maximal independent
sets in an associated graph. We first describe these graphs.

\begin{defn}
Let $D\subseteq\ZZ^d$ be finite and let $S\subseteq\ZZ^d$ satisfy $S=-S$.
The \emph{truncated Cayley graph} $\Cay(D,S)$ is the simple graph on vertex
set $D$ in which distinct $x,y\in D$ are adjacent if $x-y\in S$.
\end{defn}
Throughout the paper, all graphs are simple and self-loops are omitted.
In particular, $\Cay(D,Z-Z)=\Cay(D,(Z-Z)\setminus\{0\})$.
The vertex set $D$ will usually be $[n]^d$.

\begin{lemma}\label{lem:cover-mis}
Let $A\subseteq D\subseteq\ZZ^d$ and $Z\subseteq\ZZ^d$ be finite, with $0\in Z$.
If $D\subseteq A\oplus Z$, then $A$ is a maximal independent set in
$\Cay(D,Z-Z)$.
\end{lemma}
\begin{proof}
If distinct $a,a'\in A$ were adjacent, we would have $a-a'=z-z'$ for some
$z,z'\in Z$, and hence $a+z'=a'+z$, contradicting the disjointness property
of tilings. Thus $A$ is independent. Every $x\in D\setminus A$ has a
representation $x=a+z$, with $a\in A$ and $z\in Z\setminus\{0\}$, as
$D\subseteq A\oplus Z$. Since $x-a=z-0\in Z-Z$, the vertex $x$
has a neighbour in $A$. Therefore $A$ is maximal.
\end{proof}

Write $\mis(G)$ for the number of maximal independent sets of $G$.
We recall the following theorem on maximal independent sets in graphs.
\begin{theorem}[Moon--Moser~\cite{MoonMoser1965}]\label{thm:moon-moser}
Every graph on $N$ vertices has at most $3^{N/3}$ maximal independent sets.
\end{theorem}

Before proving the upper bound, we record that small sets contribute negligibly.
We will also use it in Sections~\ref{sec:upper-one} and~\ref{sec:upper-two}.
\begin{obs}\label{obs:small-sets}
For every fixed $d\ge1$, the number of subsets $A\subseteq[n]^d$ with
$|A|\le n^d/10$ is $o(3^{n^d/3})$.
\end{obs}
\begin{proof}
Put $N=n^d$ and let $H_2(x)=-x\log_2x-(1-x)\log_2(1-x)$.
The binary entropy bound gives
\begin{equation}\label{eq:small-sets}
    \sum_{j\le N/10}\binom Nj
    \le 2^{H_2(1/10)N}=o(3^{N/3}),
\end{equation}
since $H_2(1/10)<\frac13\log_2 3$.
\end{proof}

\begin{proof}[Proof of Theorem~\ref{thm:general-upper}]
By Observation~\ref{obs:small-sets}, it suffices to count tiles
$A\subseteq[n]^d$ with $|A|>n^d/10$. Choose $Z$ with $\ZZ^d=A\oplus Z$.
Replacing $Z$ by $Z-z$ for some $z\in Z$, we may assume that $0\in Z$.
Let $Z'\subseteq Z$ consist of those $z$ for which $(A+z)\cap[n]^d\ne\emptyset$.
Then $0\in Z'$, $[n]^d\subseteq A\oplus Z'$, and
$Z'\subseteq[-n+1,n-1]^d$. The disjoint translates $A+z$, $z\in Z'$, lie
in a box of size at most $(3n)^d$, so
\[
    |A||Z'|\le(3n)^d,\qquad |Z'|\le10\cdot3^d.
\]
Thus there are $n^{O_d(1)}$ possible sets $Z'$.

By Lemma~\ref{lem:cover-mis}, every tile counted for a given $Z'$ is a maximal
independent set in $\Cay([n]^d,Z'-Z')$. The Moon--Moser bound implies that,
for a fixed $Z'$, there are at most $3^{n^d/3}$ tiles with this specific $Z'$.
Summing over $Z'$ proves the upper bound.
\end{proof}

\subsection{Stability and tighter upper bounds}

To obtain tighter upper bounds in dimensions $1$ and $2$, we must revisit
the counting argument in Theorem~\ref{thm:general-upper}. We will need to identify the translation sets $Z'$ that contribute nearly $3^{n^d/3}$ tiles and show that there are a few choices for these translation sets.

The starting point for such an identification will be a stability result for the Moon--Moser theorem. Moon and Moser showed that, when $3\mid N$, equality
in their bound holds precisely for graphs whose connected components
are all triangles (they also established the equality cases for their
sharper bounds when $N\equiv1,2\pmod3$). We call a graph whose connected
components are all triangles a \emph{perfect triangle matching}.

\begin{defn}
An \emph{induced triangle matching} in $G$ is an induced subgraph that is
a perfect triangle matching. The \emph{triangle defect} of $G$ is
\[
    \triangledef(G)\defeq |V(G)|-
    \max\{|U|:G[U]\text{ is a perfect triangle matching}\}.
\]
\end{defn}

Following a question of Y. Rabinovich, Kahn and Park proved a stability version of the Moon--Moser theorem. This was subsequently improved quantitatively by Palmer and Patk{\'o}s~\cite{PalmerPatkos2023}, but the Kahn-Park theorem suffices for our purposes. It states that if a graph has a large number of maximal independent sets, then it contains a large induced subgraph that is a perfect triangle matching.
More precisely:

\begin{theorem}[Kahn--Park~\cite{KahnPark2020StabilityMIS}]
\label{thm:kahn-park}
There is an absolute constant $c>0$ such that every $N$-vertex graph $G$
satisfies
\[
    \mis(G)\le \exp(-c\triangledef(G))3^{N/3}.
\]
\end{theorem}

Thus, a translation set $Z'$ can contribute nearly $3^{n^d/3}$ tiles
only if its truncated Cayley graph contains an induced triangle matching
covering nearly all its vertices. We will therefore study which
translation sets can produce this structure. This will require a
detailed geometric analysis of the translation vectors and the triangles
they form.

For graphs of bounded degree, we can also ensure that these triangles
are isolated, at the cost of increasing the number of missed vertices
by a constant factor. We call $U\subseteq V(G)$ \emph{isolated} if no edge
joins $U$ to $V(G)\setminus U$.

\begin{cor}\label{cor:isolation}
If $\Delta(G)\le s$ and $\triangledef(G)=m$, then there is an isolated set
$U\subseteq V(G)$ such that $G[U]$ is a perfect triangle matching and
$|V(G)\setminus U|\le(1+3s)m$.
\end{cor}
\begin{proof}
Choose an induced triangle matching missing $m$ vertices. At most $s m$
of its triangles have an edge to a missed vertex. Delete those triangles.
The remaining triangles are isolated, and the deleted vertices together
with the originally missed vertices number at most $(1+3s)m$.
\end{proof}

\section{Typical structure in dimension one}\label{sec:upper-one}
We prove Theorem~\ref{thm:one-dimensional}, determining the number of tiles
contained in $[n]$ up to an additive $o(3^{n/3})$ error. The main contribution
comes from tiles that tile $\ZZ$ by an arithmetic progression. We first count
these tiles, then show that non-AP tiles form a negligible family.

\begin{prop}\label{prop:ap-count}
Let $r_k=n-k\lfloor n/k\rfloor$. Then
\[
    T_{\mathrm{AP}}(n,1)=
    \sum_{k=1}^n\left\lfloor\frac nk\right\rfloor^{k-r_k}
       \left(\left\lfloor\frac nk\right\rfloor+1\right)^{r_k}.
\]
Moreover, if $n=3m+s$ with $s\in\{0,1,2\}$, then
\[
    T_{\mathrm{AP}}(n,1)=(c_s+o(1))3^m,
    \qquad
    (c_0,c_1,c_2)=\left(\frac{217}{17},\frac{312}{17},\frac{450}{17}\right).
\]
\end{prop}
\begin{proof}
A set tiles by $k\ZZ$ precisely when it contains exactly one representative
of each residue class modulo $k$. Of these classes, $r_k$ have
$\lfloor n/k\rfloor+1$ representatives in $[n]$ and the rest have
$\lfloor n/k\rfloor$. Multiplying their sizes gives the summand. The families
for different $k$ are disjoint, since every set in the family for $k$ has
size $k$. This proves the exact formula.

For the asymptotics, set $F_n(k)=q^{k-r}(q+1)^r$, where
$q=\lfloor n/k\rfloor$ and $r=n-qk$. If $q=1$, then $k>n/2$ and
$F_n(k)=2^{n-k}\le2^{n/2}$. If $q\ge4$, AM--GM and the monotonicity of
$x^{1/x}$ for $x\ge3$ give $F_n(k)\le(n/k)^k\le4^{n/4}$.
Thus these two ranges contribute $O(n2^{n/2})=o(3^m)$ many tiles.

For $q=3$, the range $n/4<k\le n/3$ can be written as $k=m-j$,
where $0\le j\le J_n=m-\lfloor n/4\rfloor-1$. Since $r=s+3j$,
\[
    F_n(m-j)=3^{m-s-4j}4^{s+3j}
    =a_s3^m(64/81)^j,
    \qquad (a_0,a_1,a_2)=(1,4/3,16/9).
\]
As $J_n\to\infty$, summing this geometric progression gives
$(81a_s/17+o(1))3^m$.

For $q=2$, write $k=m+1+j$, where
$0\le j\le J'_n=\lfloor n/2\rfloor-m-1$. Here $r=m+s-2-2j$, so
\[
    F_n(m+1+j)=2^{3-s+3j}3^{m+s-2-2j}
    =b_s3^m(8/9)^j,
    \qquad (b_0,b_1,b_2)=(8/9,4/3,2).
\]
Since $J'_n\to\infty$, this range contributes $(9b_s+o(1))3^m$.
Adding the contributions gives $c_s=81a_s/17+9b_s$, as required.
\end{proof}

\begin{lemma}\label{lem:generator-stability}
For all $\rho\in (0,0.01)$, there exists $\eta=\eta(\rho)>0$ and $n_0=n_0(\rho)\in\mathbb{N}$ such that the
following holds for $n\geq n_0$. If $R\subseteq[n-1]$ and
\[
    \triangledef\bigl(\Cay([n],R\cup(-R))\bigr)<\eta n,
\]
then there is an integer $k\ge1$ such that
\begin{equation}\label{eq:one-dimensional-generators}
    |k-n/3|\le\rho n,
    \qquad R\cap[1,(1-\rho)n]=\{k,2k\}.
\end{equation}
\end{lemma}
\begin{proof}
Fix $\rho\in (0,0.01)$ and let $\eta=\rho/10$.
Let $G=\Cay([n],R\cup(-R))$, and let $U$ be the vertices missed by a
largest induced triangle matching, so $|U|<\eta n$.
If $R$ is empty, then $|U|=n$, a contradiction. Enumerate
$R=\{x_1<x_2<\cdots\}$.

Suppose $x_3$ exists and $x_3\le(1-\rho)n$. For each integer
$1\le t\le n-x_3$, the vertices $t,t+x_1,t+x_2,t+x_3$ cannot all lie
outside $U$, since $t$ would have degree at least three in the matching.
Each vertex occurs in at most four such quadruples, giving
$|U|\ge(n-x_3)/4\ge\rho n/4$, a contradiction. Hence, if present,
$x_3>(1-\rho)n$.

If $x_1<n/3-\rho n$, then for each integer
$1\le t\le\lfloor\rho n\rfloor$, the four vertices
$t,t+x_1,t+2x_1,t+3x_1$ lie in $[n]$ and induce a connected graph.
They cannot all belong to an induced triangle matching. The same counting argument
gives $|U|\ge\lfloor\rho n\rfloor/4>\eta n$ for sufficiently large $n$.
If $x_1>n/3+\rho n$, every integer
$t\in[n/3-\rho n,n/3+\rho n]$ lies in no triangle. Indeed, $t<x_1$,
so it must be the smallest vertex of any such triangle. The other vertices
would be $t+r,t+s$, with $r\ge x_1$ and $s-r\ge x_1$, forcing
$t+s\ge t+2x_1>n$. This again gives $|U|>\eta n$.
Thus $|x_1-n/3|\le\rho n$.

Finally, suppose $|R|=1$ or $x_2\ne2x_1$. Every integer
$t\in[2\rho n,3\rho n]$ satisfies $t<x_1$ and cannot use a positive
edge length exceeding $(1-\rho)n$. A triangle containing $t$ would therefore
require $x_2$ to exist and $x_2-x_1\in R$, forcing $x_2=2x_1$,
a contradiction. There are at least $\rho n-1$ such
integers, so $|U|>\eta n$. We conclude that $x_2=2x_1$.
Taking $k=x_1$ proves~\eqref{eq:one-dimensional-generators}, since
$2k\le(2/3+2\rho)n<(1-\rho)n$.
\end{proof}

The second lemma uses the disjoint covering of a long interval to turn the
generator condition into an AP tiling.
\begin{lemma}\label{lem:generators-force-ap}
There exist absolute constants $\rho>0$ and $L,n_0\in\mathbb N$ such that the
following holds for all $n\geq n_0$. Suppose $A\subseteq[n]$ and
$Z'\subseteq\ZZ$ is finite, with $0\in Z'$ and
$[Ln]\subseteq A\oplus Z'$. If $R=(Z'-Z')\cap[n-1]$ satisfies
\eqref{eq:one-dimensional-generators} for an integer $k\ge1$, then $A$
contains exactly one representative of each residue class modulo $k$.
Consequently, $\ZZ=A\oplus k\ZZ$.
\end{lemma}
\begin{proof}
Choose $\rho>0$ sufficiently small, and put
$\alpha=1/3+2\rho$ and $\beta=1/3-\rho$. Since $4\beta>3\alpha$, we can
fix an integer $L$ such that $\frac{4\beta}{3}(L/\alpha-3)>L+1$.
By Lemma~\ref{lem:cover-mis}, $A$ is a maximal independent set in
$G=\Cay([n],Z'-Z')$.

Since $k>n/4$, each residue class of $[n]$ modulo $k$ has at most four
elements. Independence and $k,2k\in R$ imply that $A$ meets each class at
most once, except possibly the $(n-3k)_+$ classes of size four, which it
meets at most twice. Hence
\[
    |A|\le k+(n-3k)_+\le\alpha n.
\]
Enumerate the translates meeting $[Ln]$ as
$Y=\{z\in Z':(A+z)\cap[Ln]\ne\emptyset\}=\{y_1<\cdots<y_m\}$.
Note that $Y\oplus A\supseteq[Ln]$ and $1-n\le y_i\le Ln-1$, so
$m|A|\ge Ln$ and $y_m-y_1\le(L+1)n$. In particular, $m\ge L/\alpha$.
By~\eqref{eq:one-dimensional-generators}, every successive gap is either
$k$ or at least $2k$. If no three consecutive gaps equalled $k$, every
block of three gaps would have total length at least $4k$. Therefore
\[
    (L+1)n\ge y_m-y_1
    \ge4k\left\lfloor\frac{m-1}{3}\right\rfloor
    \ge\frac{4\beta}{3}(L/\alpha-3)n,
\]
contrary to the choice of $L$. Thus $z,z+k,z+2k,z+3k\in Z'$ for some $z$.
Their differences give $k,2k,3k\in Z'-Z'$. Since $k>n/4$, every residue
class of $[n]$ modulo $k$ is now a clique in $G$, and $A$ meets it at most
once.

Let $C$ be a residue class modulo $k$, and take its unique element
$x\in[k+1,2k]\subseteq[n]$. Since
$\max\{x-1,n-x\}\le(2/3+2\rho)n<(1-\rho)n$,
condition~\eqref{eq:one-dimensional-generators} implies that every neighbour
of $x$ is obtained using $k$ or $2k$, and hence belongs to $C$.
If $A$ missed $C$, then $x\notin A$ would have no neighbour in $A$,
contradicting maximality. Thus $A$ meets every residue class exactly once,
as required.
\end{proof}

\begin{proof}[Proof of Theorem~\ref{thm:one-dimensional}]
Fix $\rho,L$ as in Lemma~\ref{lem:generators-force-ap}, and let $\eta$ be
given by Lemma~\ref{lem:generator-stability} for this $\rho$.
Observation~\ref{obs:small-sets} bounds the number of tiles of size at most $n/10$.
For a remaining tile $A$, choose $Z$ with $\ZZ=A\oplus Z$ and $0\in Z$,
and let $Z'\subseteq Z$ consist of the translates meeting $[Ln]$.
Then $0\in Z'$, $Z'\subseteq[1-n,Ln-1]$, and $[Ln]\subseteq A\oplus Z'$.
These disjoint translates lie in an interval of size at most $(L+2)n$, so
$|A||Z'|\le(L+2)n$. Since $|A|>n/10$, we have $|Z'|\le10(L+2)$.
There are therefore $n^{O(L)}$ possible sets $Z'$.

For a fixed $Z'$, put $G=\Cay([n],Z'-Z')$ and
$R=(Z'-Z')\cap[n-1]$, so $G=\Cay([n],R\cup(-R))$.
By Lemma~\ref{lem:cover-mis}, every corresponding tile is a maximal
independent set in $G$. If $\triangledef(G)<\eta n$,
Lemma~\ref{lem:generator-stability} gives~\eqref{eq:one-dimensional-generators},
and Lemma~\ref{lem:generators-force-ap} forces an AP tiling. Thus any $Z'$
arising from a non-AP tile has $\triangledef(G)\ge\eta n$.
Applying Theorem~\ref{thm:kahn-park} and summing over $Z'$ gives
\[
    T_{\mathrm{non\text{-}AP}}(n,1)
    \le\sum_{j=0}^{\lfloor n/10\rfloor}\binom nj
       +n^{O(L)}e^{-c\eta n}3^{n/3}
    =o(3^{n/3}).
\]
Combining this with Proposition~\ref{prop:ap-count}, for $n=3m+s$ we obtain
$T(n,1)=T_{\mathrm{AP}}(n,1)+o(3^{n/3})=(c_s+o(1))3^m$, as required.
\end{proof}
\section{Near-tight bounds in dimension two}\label{sec:upper-two}

We prove Theorem~\ref{thm:two-dimensional} in this section. Using Lemma~\ref{lem:cover-mis} and Theorem~\ref{thm:kahn-park}, the proof reduces to showing that there are at most $\widetilde{O}(n)$ choices for $S\subset [n]^2$ of bounded size satisfying $\triangledef(G)=O(\log n)$, where $G=\Cay([n]^2,S\cup -S)$. We will divide the generators into heavy generators and non-heavy generators, depending on if they meet more than $(\log n)^4$ vertices in $G$. There are $\widetilde{O}(1)$ choices for non-heavy generators. On the other hand, heavy generators must be part of a triangle whose translates meet at least $(\log n)^3$ vertices of $G$. We will call the family of such triangles $\mathcal{T}_{\mathrm{heavy}}$ \footnote{the precise definition of $\mathcal{T}_{\mathrm{heavy}}$ involves triangle types instead of triangles, which is defined later in the section}. Most of the work in this section will go into identifying the geometric structure of $\mathcal{T}_{\mathrm{heavy}}$ enough to conclude that there are only $\widetilde O(n)$ possibilities for it, which will conclude the proof.

\subsection{Triangle types and their supports}

Throughout this section, $[a]=\emptyset$ for integers $a\le0$.

A \emph{triangle type} is a multiset $T=\{a,b,c\}$ of nonzero vectors in
$\ZZ^2$ satisfying $a+b+c=0$, where we identify $T$ and $-T\defeq \{-a,-b,-c\}$.
Thus a triple $\triangle=\{x,y,z\}\subseteq[n]^2$ determines its type $T=T(\triangle)\defeq \{x-y,y-z,z-x\}$. For an ordering $(a,b,c)$ of a representative of $T$, define
\[
 R^n_{(a,b,c)}=
 \bigcup_{\substack{x\in[n]^2:\\x+a,x+a+b\in[n]^2}}
       \{x,x+a,x+a+b\},
 \qquad
 S_T^n=R^n_{(a,b,c)}\cup R^n_{(-a,-b,-c)}.
\]
We call $S_T^n$ the \emph{support} of $T$. Note that the support is well-defined as $R^n_{(\sigma(a),\sigma(b),\sigma(c))}=R^n_{\operatorname{sgn}(\sigma)(a,b,c)}$ for $\sigma\in \text{Sym}(\{a,b,c\})$. It  is precisely the union of all triples of type $T$ contained in $[n]^2$. 

Define the \emph{crossants}
\[
 C_+=\{(x,y):xy\ge0\},\qquad C_-=\{(x,y):xy\le0\}.
\]
A triangle type is of \emph{$D$-type} if its vectors all belong to $C_+$ or all
belong to $C_-$, and of \emph{$L$-type} otherwise. 
Write
\[
 \rho(x_1,x_2)=(n+1-x_1,n+1-x_2)
\]
for central reflection of the square $[n]^2$.
Similarly, write
\[
r_1(x_1,x_2)=(n+1-x_1,x_2)\qquad\text{and}\qquad r_2(x_1,x_2)=(x_1,n+1-x_2)
\]
for the two \emph{coordinate reflections} of the square. For an \emph{edge vector} $u=x-y=(u_1,u_2)\in \mathbb{Z}^2$ with $x,y\in [n]^2$, we get the induced central reflection and coordinate reflections as
\[
\widetilde{\rho}(u)=-u, \widetilde{r}_1(u)=(-u_1,u_2)\text{ and } \widetilde{r}_2(u)=(u_1,-u_2).
\]

\begin{lemma}\label{lem:two-supports}
After reflecting a coordinate (that is, after applying $\widetilde{r}_1$ to all vectors in the tuple) if necessary, every triangle type has one
of the following representations.
\begin{enumerate}
\item If it is of $D$-type, then
$T=\{u,v,-u-v\}$, with nonzero $u,v\in\ZZ_{\ge0}^2$, and
\begin{equation}\label{eq:D-boxes}
 S_T^n=B\cup(B+u)\cup(B+v)\cup(B+u+v),\qquad
 B=[n-u_1-v_1]\times[n-u_2-v_2].
\end{equation}
\item If it is of $L$-type, then
$T=\{u,-v,v-u\}$, where
$u=(P,y)$, $v=(x,Q)$, $P>x\ge0$ and $Q>y\ge0$, and
\begin{equation}\label{eq:L-boxes}
 S_T^n=R\cup\rho(R),\qquad
 R=B\cup(B+u)\cup(B+v),\qquad B=[n-P]\times[n-Q].
\end{equation}
\end{enumerate}
\end{lemma}
\begin{proof}
For a $D$-type, reflect a coordinate so that all three vectors lie in
$C_+$. Each nonzero vector in $C_+$ lies in exactly one of the nonnegative
and nonpositive quadrants. Two of the three lie in the same quadrant;
changing all signs if necessary, call these $u,v\in\ZZ_{\ge0}^2$.
The third is $-u-v$. The two triangle families are
\[
 \{z,z+u,z+u+v\},\qquad \{z,z+v,z+u+v\},\qquad z\in B.
\]
Indeed, both triangles lie entirely in $[n]^2$ precisely when
$z\in B=[n-u_1-v_1]\times[n-u_2-v_2]$. This proves~\eqref{eq:D-boxes}.

For an $L$-type, two of the three vectors lie in the same crossant,
since $C_+\cup C_-=\ZZ^2$. Reflect a coordinate so that these two
vectors lie in $C_+$. They cannot lie in the same quadrant, since their
negative sum would then lie in the opposite quadrant, making all three
vectors lie in $C_+$ and $T$ of $D$-type. Thus, relabelling these two
vectors, we can write them as $u,-v$, with $u,v\in\ZZ_{\ge0}^2$.
The third is $v-u$, which must lie outside $C_+$, so
$(u_1-v_1)(u_2-v_2)<0$. Interchanging $u,v$ and changing all signs if necessary, we obtain $u_1>v_1$ and $u_2<v_2$. Write
$u=(P,y)$ and $v=(x,Q)$, so $P>x\ge0$ and $Q>y\ge0$.
The two triangle families are
\[
 \{z,z+u,z+v\},\qquad
 \{\rho(z),\rho(z+u),\rho(z+v)\},\qquad z\in B.
\]
The first triangle lies entirely in $[n]^2$ precisely when
$z\in B=[n-P]\times[n-Q]$. Central reflection gives exactly the
family with all edge signs reversed,  proving
\eqref{eq:L-boxes}.
\end{proof}

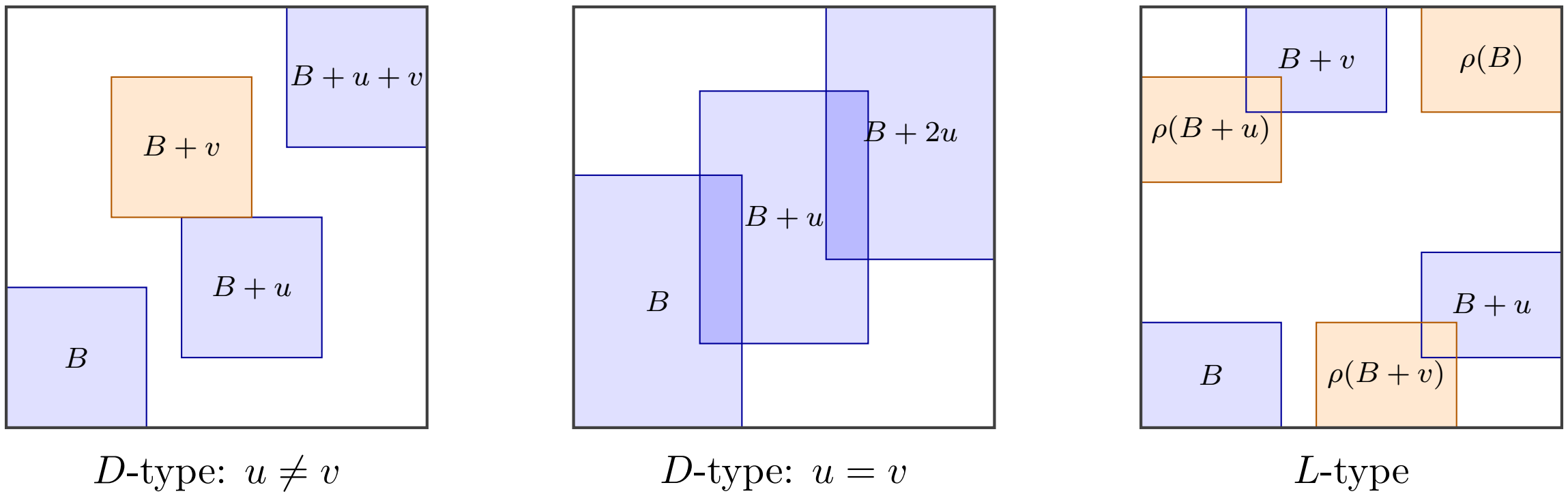
\begin{figure}[t]
\centering
\begin{tikzpicture}[x=0.34cm,y=0.34cm,every node/.style={font=\scriptsize}]
 \begin{scope}
  \filldraw[fill=blue!12,draw=blue!60!black] (0,0) rectangle (4,4);
  \filldraw[fill=blue!12,draw=blue!60!black] (5,2) rectangle (9,6);
  \filldraw[fill=orange!18,draw=orange!70!black] (3,6) rectangle (7,10);
  \filldraw[fill=blue!12,draw=blue!60!black] (8,8) rectangle (12,12);
  \draw[black!75,thick] (0,0) rectangle (12,12);
  \node at (2,2) {$B$};
  \node at (7,4) {$B+u$};
  \node at (5,8) {$B+v$};
  \node at (10,10) {$B+u+v$};
  \node[font=\small] at (6,-1.3) {$D$-type: $u\neq v$};
 \end{scope}
 \begin{scope}[xshift=5.5cm]
  \fill[blue!12] (0,0) rectangle (4.8,7.2);
  \fill[blue!12] (3.6,2.4) rectangle (8.4,9.6);
  \fill[blue!12] (7.2,4.8) rectangle (12,12);
  \fill[blue!27] (3.6,2.4) rectangle (4.8,7.2);
  \fill[blue!27] (7.2,4.8) rectangle (8.4,9.6);
  \draw[blue!60!black] (0,0) rectangle (4.8,7.2);
  \draw[blue!60!black] (3.6,2.4) rectangle (8.4,9.6);
  \draw[blue!60!black] (7.2,4.8) rectangle (12,12);
  \draw[black!75,thick] (0,0) rectangle (12,12);
  \node at (2.4,3.6) {$B$};
  \node at (6,6) {$B+u$};
  \node at (9.6,8.4) {$B+2u$};
  \node[font=\small] at (6,-1.3) {$D$-type: $u=v$};
 \end{scope}
 \begin{scope}[xshift=11cm]
  \fill[blue!12] (0,0) rectangle (4,3);
  \fill[blue!12] (8,2) rectangle (12,5);
  \fill[blue!12] (3,9) rectangle (7,12);
  \fill[orange!18] (8,9) rectangle (12,12);
  \fill[orange!18] (0,7) rectangle (4,10);
  \fill[orange!18] (5,0) rectangle (9,3);
  \draw[blue!60!black] (0,0) rectangle (4,3);
  \draw[blue!60!black] (8,2) rectangle (12,5);
  \draw[blue!60!black] (3,9) rectangle (7,12);
  \draw[orange!70!black] (8,9) rectangle (12,12);
  \draw[orange!70!black] (0,7) rectangle (4,10);
  \draw[orange!70!black] (5,0) rectangle (9,3);
  \draw[black!75,thick] (0,0) rectangle (12,12);
  \node at (2,1.5) {$B$};
  \node at (10,3.5) {$B+u$};
  \node at (5,10.5) {$B+v$};
  \node at (10,10.5) {$\rho(B)$};
  \node at (2,8.5) {$\rho(B+u)$};
  \node at (7,1.5) {$\rho(B+v)$};
  \node[font=\small] at (6,-1.3) {$L$-type};
 \end{scope}
\end{tikzpicture}
\caption{Various possibilities for supports in positive orientation. Reflecting either coordinate gives the
negative orientation.}
\label{fig:triangle-shapes}
\end{figure}

We call $B$ the \emph{fundamental rectangle}. The displayed representations
have corner rectangles $B$ and $\rho(B)$ at the bottom-left and top-right
corners; we call this a \emph{positive orientation}. Reflecting one
coordinate gives a \emph{negative orientation}, whose corner rectangles
are at the other two corners. An axis-parallel edge can allow more than
one such representation; choose one whenever a representation is needed.
The arguments below use the actual corner rectangles, and hence apply
to every such choice. The various possibilities for the support are illustrated in Figure~\ref{fig:triangle-shapes}.

We say that $(A_i)_{i\in I}$ is an \emph{$r$-approximate partition} of
$[n]^2$ if $|A_i\cap A_j|\le r$ for distinct $i,j$ and
$|[n]^2\setminus\bigcup_i A_i|\le r$.

\begin{lemma}\label{lem:triangle-partition}
Let $G=\Cay([n]^2,S)$, where $S=-S$, $0\notin S$ and $|S|\le s$.
There is a collection $\mathcal T$ of at most $s^3$ triangle types from
$S$ and a set $X$ with
\[
 |X|\le(1+3s)\triangledef(G)
\]
such that their supports form an $|X|$-approximate partition.
Every vertex outside $X$ belongs to a unique triangle of $G$, which is
an isolated component and has type in $\mathcal T$.
\end{lemma}
\begin{proof}
The maximum degree of $G$ is at most $s$. Apply Corollary~\ref{cor:isolation}
and let $X$ be the complement of the isolated triangle matching obtained
there. Take one type, with $T$ and $-T$ identified, for each distinct type
of its components. There are at most $s^3$ choices of three edge vectors.
Their supports cover $[n]^2\setminus X$. If a vertex outside $X$ belongs
to the support of a type, the corresponding triangle must be its isolated
component. Thus it cannot belong to supports of two distinct types,
proving the intersection bound and all the assertions.
\end{proof}

There are two further consequences that we will use repeatedly. First,
if a $D$-type in $\mathcal T$ has $u\ne v$ in~\eqref{eq:D-boxes}, every
$z\in B$ has the three distinct neighbours $z+u,z+v,z+u+v$. It therefore
belongs to $X$, and
\begin{equation}\label{eq:non-ap-D}
 |S_T^n|\le4|B|\le4|X|.
\end{equation}
For a $D$-type with $u=v$, the support consists of the three rectangles
$B,B+u,B+2u$.

Second, any two of these three rectangles, or of the six rectangles for
an $L$-type, intersect only in $X$. For the three $D$-rectangles, a point
in two rectangles belongs to two distinct translates of
$\{0,u,2u\}$. For an $L$-type, two boxes in the same family similarly
give distinct translates of $\{0,u,v\}$. The two families cannot give
the same triangle: that would make a noncollinear three-point set
centrally symmetric. Central symmetry of a three-point set fixes one
point and interchanges the other two, forcing collinearity, whereas
$PQ-xy>0$. Thus in every case such a point belongs to two distinct
triangles and must lie in $X$. Rectangles from distinct types also
intersect only in $X$, by Lemma~\ref{lem:triangle-partition}.

\subsection{Heavy triangle types}
Fix constants $s,K$ and a set $S$ as in Lemma~\ref{lem:triangle-partition}
with $\triangledef(G)\le K\log n$. Lemma~\ref{lem:triangle-partition} then supplies a collection of triangle types $\mathcal{T}=\mathcal{T}(S)$.  In this subsection all implicit constants
may depend on $s,K$, and $n$ is sufficiently large. Put
$r=1+\lceil(1+3s)K\log n\rceil$, so $|X|\le r$, and set
\[
    \lambda=(\log n)^3,\qquad
    \mathcal T_{\mathrm{heavy}}
       =\{T\in\mathcal T:|S_T^n|\ge\lambda\}.
\]
Note that $\mathcal T_{\mathrm{heavy}}=\mathcal T_{\mathrm{heavy}}(S)$, but we will suppress the dependence in the notation to avoid clutter. If $T\in \mathcal T_{\mathrm{heavy}}$, then we say that the type $T$ and the support $S^n_T$ are \emph{heavy}. Since $|\mathcal T|=O(1)$, the heavy supports cover all but $O(\lambda)$
points. We repeatedly use $r^2=o(\lambda)$ and $\lambda=o(n)$.
By~\eqref{eq:non-ap-D}, every heavy $D$-type has the form
$\{u,u,-2u\}$, after reflection. Its support has three boxes; a heavy
$L$-support has six. All these boxes have pairwise intersections of
size at most $r$, so their areas satisfy
\begin{equation}\label{eq:support-areas}
 |S_T^n|=3ab+O(r)\quad(T \text{ of type }D),\qquad
 |S_T^n|=6ab+O(r)\quad(T \text{ of type } L),
\end{equation}
where $a,b$ are the dimensions of the fundamental rectangle. In particular, a heavy fundamental
rectangle has area at least $\lambda/7$ (for $n$ sufficiently large). Here and below \emph{area}
means the number of lattice points.

For two rectangles with a common corner or centre, the intersection has
side lengths equal to the smaller of the corresponding side lengths.
Thus, if one rectangle has both dimensions longer than $r$ and the other has area greater than $r^2$, their intersection has more than $r$ points.

For a positive $L$-type, the bottom and left pairs of boxes (by which we mean the pairs $(B,\rho(B+v))$ and $(B,\rho(B+u))$ respectively in Figure~\ref{fig:triangle-shapes}), and the pair $(\rho(B+v),B+u)$ give respectively
\begin{equation}\label{eq:L-overlaps}
    b(2a+x-n)_+\le r,\qquad a(2b+y-n)_+\le r,\qquad
    (a-x)_+(b-y)_+\le r.
\end{equation}
In particular, if $a,b>r$, integrality gives $2a+x\le n$ and $2b+y\le n$. The last inequality in \eqref{eq:L-overlaps} then implies

\begin{equation}\label{eq:L-overlaps2}
   (3a-n)_+(3b-n)_+\le r, \quad \text{if } a,b>r.
\end{equation}

\begin{lemma}\label{prop:heavy-structure}
Let $S\subset \mathbb{Z}^2\setminus\{0\}$ be such that $S=-S$, $|S|=O(1)$ and $G=\Cay([n]^2,S)$ satisfies $\triangledef(G)=O(\log n)$. Then the following holds for $n$ sufficiently large. There are at most two heavy types.
If there are two, one is $D$-type, one is $L$-type, and their orientations are opposite.
\end{lemma}
\begin{proof}

Some heavy support has size $\Omega(n^2)$, since there are only $O(1)$ types
and their supports cover all but $O(\lambda)$ points. Reflect the square so
that this support is positive, and let $a,b=\Omega(n)$ be its fundamental
dimensions. Every other heavy type must be negative: if not, its positive corner
rectangle would have area at least $\lambda/7>r^2$ and would intersect the
large positive corner rectangle in more than $r$ points.
If another heavy type has both fundamental dimensions greater than $r$,
its negative corner rectangle similarly excludes any third heavy type.
Otherwise, all the other supports together have size $O(nr)$, so the large
support has size $n^2-O(nr)$.

Suppose first that this large support is of $D$-type. Then, by \eqref{eq:support-areas},
$ab=n^2/3+O(nr)$, while the overlap of its first two boxes gives
$((3a-n)/2)_+((3b-n)/2)_+\le r$.
Hence $\min\{a,b\}\le n/3+O(\sqrt r)$. Combining this with
$\min\{a,b\}\ge ab/n$, and interchanging coordinates if necessary, gives
$a=n/3+O(r)$ and $b=n-O(r)$.

If the large support is of $L$-type, then $ab=n^2/6+O(nr)$.
By~\eqref{eq:L-overlaps2},
$(3a-n)_+(3b-n)_+\le r$, so again $\min\{a,b\}\le n/3+O(\sqrt r)$.
Take $b$ to be the larger dimension. The area bound gives
$b\ge n/2-O(r)$, and $2b+y\le n$ then gives
$b=n/2+O(r)$, $y=O(r)$ and $a=n/3+O(r)$.
Since $b-y>r$, the last inequality in~\eqref{eq:L-overlaps} gives $x\ge a$.
Together with $2a+x\le n$, this yields $x=n/3+O(r)$.

In either case, the large support contains a rectangle along the left side
of height $\Omega(n)$ whose top edge is at height $n-e$, where $e=O(r)$.
(For $D$, this is $B$, whereas for $L$, this is $\rho(B+u)$ in Figure~\ref{fig:triangle-shapes}.)
Any other heavy type has a top-left corner rectangle, say of dimensions
$c,d$. Its intersection with this left-side rectangle gives $d\le e+r=O(r)$.
Since $cd\ge\lambda/7$, its width satisfies $c>r$.
Two such corner rectangles would intersect in more than $r$ points.
This proves that there are at most two heavy types, with opposite orientations.
Two $D$-types are impossible as well: their middle boxes have a common centre,
one has both dimensions $\Omega(n)$, and the other has area at least
$\lambda/7>r^2$.

It remains to exclude two $L$-types. Let their fundamental dimensions be
$(a,b)$ and $(c,d)$, with the first in positive orientation. If all four dimensions were
smaller than $n/2-C\sqrt\lambda$, then $[n/2-0.1C\sqrt\lambda, n/2+0.1C\sqrt\lambda]^2$ would be disjoint from both heavy supports, thus exceeding the $O(\lambda)$ uncovered allowance for sufficiently large $C$.
It follows that one of the dimensions exceeds $n/2-O(\sqrt\lambda)$, say, $b\ge n/2-O(\sqrt\lambda)$.
The second bound in~\eqref{eq:L-overlaps} gives
$b\le n/2+O(r)$ and $y=O(\sqrt\lambda)$.
Its upper-left box ends at height $n-y$, so the negative corner rectangle
has $d\le y+r=O(\sqrt\lambda)$. Hence $c\ge\lambda/(7d)> r$.
The negative support has size $O(n\sqrt\lambda)$, giving
$ab=n^2/6+O(n\sqrt\lambda)$ and $a=n/3+O(\sqrt\lambda)$. Now, considering how $[n]\times \{1\}$ is covered by the two supports, we obtain that $2a+2c=n+O(\lambda)$, so $c=n/6+O(\lambda)$. Since $\min\{a,c\}>r$, the bottom-right corner rectangle of the negative support cannot overlap with the positive
lower-right rectangle even in one row, so $y\ge d\ge 1$.

Consider the set of points $J=[0.4n,0.6 n]\times [b+1,n-b]$. Note that $|J|=\Omega(n)$ as $n\geq 2b+y\geq 2b+1$. We will show that $J$ does not intersect either of the heavy supports, which would yield the required contradiction. This is because any of the four boxes of either support touching the bottom or top side of $[n]^2$ do not intersect $\mathbb{Z}\times [b+1,n-b]$. The remaining two boxes of either support touch the left or the right side of $[n]^2$ and hence do not intersect $[0.4n, 0.6n]\times\mathbb{Z}$.
\end{proof}

\begin{prop}\label{prop:heavy-count}
\[
\#\{\mathcal T_{\mathrm{heavy}}(S):S\subset\mathbb{Z}^2\setminus\{0\},S=-S, |S|\leq s, \triangledef(\Cay([n]^2,S))\leq K\log n\}=\widetilde{O}_{s,K}(n).
\]
\end{prop}
\begin{proof}

By Lemma~\ref{prop:heavy-structure}, we know that for the relevant $S$, $|\mathcal{T}_{\mathrm{heavy}}|\leq 2$.
Suppose first that $|\mathcal{T}_{\mathrm{heavy}}|=1$. With only one heavy type, its support has size $n^2-O(\lambda)$.
Since $\lambda=o(nr)$, the estimates for a support of size $n^2-O(nr)$
in the preceding proof apply. For $D$-type, they give
$a=n/3+O(r)$ and $b=n-O(r)$, after interchanging coordinates if necessary.
For $L$-type, they give $a=n/3+O(r)$, $b=n/2+O(r)$,
$x=n/3+O(r)$ and $y=O(r)$.
Therefore, after the appropriate reflection,
a single heavy $D$-type satisfies
\[
    \{n-2u_1,n-2u_2\}
      =\{n/3+O(\lambda),\ n-O(\lambda)\},
\]
and a single heavy $L$-type satisfies, after interchanging coordinates,
\[
    (n-u_1,n-v_2,u_1-v_1,v_2-u_2)
      =(n/3,n/2,n/3,n/2)+O(\lambda).
\]
In both cases, there are $O(\lambda^4)$ choices for $\mathcal{T}_{\mathrm{heavy}}$.

When $|\mathcal{T}_{\mathrm{heavy}}|= 2$, we know from Lemma~\ref{prop:heavy-structure} that one is $D$-type and the other is $L$-type with opposite orientations. Reflect the square so that the $D$-type is positive, and write
$u=(p,q)$, $v=(-P,t)$ and $w=(-x,Q)$, where $P>x\ge0$ and $Q>t\ge0$.
Their fundamental dimensions are $a=n-2p$, $b=n-2q$, $c=n-P$ and $d=n-Q$.
The vector $z=w-v=(P-x,Q-t)$ is positive in both coordinates.
It gives an edge starting at every point of
$[n-P+x]\times[n-Q+t]$, a corner rectangle containing $[c]\times[d]$.
This rectangle and the $D$-rectangle $[a]\times[b]$ both have area at least
$\lambda/7$, and one has area $\Omega(n^2)$. Their intersection therefore
contains a point outside $X$. Its edge of difference $z$ must lie in its
isolated $D$-triangle, so $z\in\{u,2u\}$.

Fix $h\in\{1,2\}$ with $w-v=hu$. The dimensions $a,b,c,d$ determine both
types, since $p=(n-a)/2$, $q=(n-b)/2$, $P=n-c$, $Q=n-d$,
$x=P-hp$ and $t=Q-hq$. We have $p,q>0$ because $z$ is positive in both
coordinates. Also $x,t>0$: if not, we would get one of the six translates of the fundamental rectangle of the negative $L$-type at
the bottom-left corner, which would intersect the fundamental rectangle of the positive $D$-type in more than $r$ points.
Now, considering how $[n]^2$, $\{1\}\times [n]$ and $[n]\times \{1\}$ are respectively covered by both supports, we obtain that
\begin{equation}\label{eq:mixed-dimensions}
    3ab+6cd=n^2+O(\lambda),\qquad
    b+2d=n+O(\lambda),\qquad
    a+2c=n+O(\lambda).
\end{equation}
Here the $D$-support contributes one box to each of the sides and the
$L$-support contributes two as $p,q,x,t,P,Q>0$. The errors account for the box intersections
and the $O(\lambda)$ uncovered points.

There are $O(\lambda^3)$ choices for the three integer right-hand sides -- fix them. If $a=c$, the last equation already fixes both, and choosing $b$ in at most $n$ ways determines $d$.
Otherwise, choose $a$ in at most $n$ ways. The last equation fixes $c$.
As $a\ne c$, the remaining two linear equations determine $b,d$ uniquely.
Thus there are $O(n\lambda^3)$ possible choices for $a,b,c,d$ and hence also for $\mathcal{T}_{\mathrm{heavy}}$,
completing the proof.
\end{proof}

\subsection{Counting tiles}
For $z\in[-n+1,n-1]^2\setminus\{0\}$, define
\[
    E_n(z)=\{(x,x+z):x,x+z\in[n]^2\},\qquad
    q_n(z)=|E_n(z)|=(n-|z_1|)(n-|z_2|).
\]
Call $z$ a \emph{heavy generator} if $q_n(z)\ge(\log n)^4$.

\begin{proof}[Proof of Theorem~\ref{thm:two-dimensional}]
Tiles of size at most $n^2/10$ are negligible in number by Observation~\ref{obs:small-sets}.
Recall from the proof of Theorem~\ref{thm:general-upper} that every remaining
tile is a maximal independent set in $G=\Cay([n]^2,Z'-Z')$ for some $Z'$ of
bounded size. Equivalently, $G=\Cay([n]^2,S)$, where
$S=((Z'-Z')\cap[-n+1,n-1]^2)\setminus\{0\}$ is symmetric and $|S|=O(1)$.
There are $n^{O(1)}$ possible choices for $S$. Theorem~\ref{thm:kahn-park} shows
that sets $S$ with $\triangledef(G)>K\log n$ contribute at most
$n^{O(1)-cK}3^{n^2/3}=o(3^{n^2/3})$ for some sufficiently large constant $K$.

For the remaining sets $S$, we apply Lemma~\ref{lem:triangle-partition} and
Proposition~\ref{prop:heavy-count}.
We claim that for $S$ satisfying $\triangledef(G)\leq K\log n$, there are $\widetilde O(1)$ choices for $S$ having fixed $\mathcal{T}_{\mathrm{heavy}}(S)$. As there are $\widetilde O(n)$ choices for relevant $\mathcal{T}_{\mathrm{heavy}}(S)$ by Proposition~\ref{prop:heavy-count}, the claim implies that there are $\widetilde O(n)$ choices for sets $S$ so that $\triangledef(G)\leq K\log n$. Having fixed such an $S$, $G$ has at most $3^{n^2/3}$
maximal independent sets by Theorem~\ref{thm:moon-moser}, which establishes
\[
T(n,2)=\widetilde{O}(n3^{n^2/3}).
\]
It therefore only remains to establish the claim.

Every heavy generator in $S$ occurs, up to sign, in a heavy triangle type. If $z$ occurs in no heavy triangle type, every edge of difference $z$ is
either inside a triangle of a nonheavy type or has an endpoint in $X$, by Lemma~\ref{lem:triangle-partition}.
The former edges number $O(\lambda)$, because there are only $O(1)$ types;
the latter number $O(|X|)$, because the degree is bounded. Thus
$q_n(z)=O(\lambda)<(\log n)^4$.

Therefore, having fixed $\mathcal{T}_{\mathrm{heavy}}$, there are only $O(1)$ choices for the heavy generators in $S$. Finally, if $z\in [-n+1,n-1]^2\setminus\{0\}$ is not heavy, then $|z_1|,|z_2|\in [n-(\log n)^4,n]$ and so there are only $\widetilde{O}(1)$ choices for a single non-heavy generator, and hence also for all the non-heavy generators in $S$ as $|S|=O(1)$.

\end{proof}

\section{Lower bound}\label{sec:lower}
We prove Theorem~\ref{thm:general-lower} by constructing $n^{\binom d2}$
lattices, each giving $3^{n^d/3}$ tiles in $[n]^d$. We show that these
families are mostly disjoint: almost every set counted belongs to exactly
one family. This gives at least $(1-o(1))n^{\binom d2}3^{n^d/3}$ tiles,
as required.

Suppose $3\mid n$. Throughout this section, for convenience, we count tiles
contained in $Q=\{0,\ldots,n-1\}^d$, which are in bijection with tiles
contained in $[n]^d$. Let $\mathcal M$ be the family of upper triangular
integer matrices with diagonal $(n,\ldots,n,n/3)$ and entries
$0\le M_{ij}<n$ for $i<j$. Thus $|\mathcal M|=n^{\binom d2}$.
Let $L_M$ be the lattice generated by the columns of $M$, and put
\[
    G_M=\Cay(Q,L_M).
\]

\begin{obs}\label{obs:lattice-triangles}
For every $M\in\mathcal M$, every coset of $L_M$ has exactly three representatives in $Q$. Thus the graph $G_M$ is a perfect triangle matching.
\end{obs}
\begin{proof}
Let $x\in \ZZ^d$ be any point, then we show that $x$ can be reduced to exactly three representatives in $Q$ by adding integer linear combinations of the columns of $M$. 

First, using the last column of $M$, we can reduce $x$ so that $0\leq x_d<n$, and there are exactly three possible choices for $x_d$ separated by $n/3$. For each choice,
reduce the remaining coordinates modulo $n$, from the $(d-1)$-th coordinate upwards to the first. The reduction on coordinate $i$ uses only column $i$ of $M$, so it leaves the
coordinates below it unchanged and determines the next coordinate uniquely.
Thus each of the three choices for $x_d$ gives exactly one point of
$Q$ in the coset.
\end{proof}

Write $\mathcal F_M$ for the family of maximal independent sets of $G_M$. By the above observation, each element $F\in \mathcal F_M$ picks one point from each triangle matching component of $G_M$, and hence one representative from each coset of $L_M$. Thus $F\oplus L_M=\ZZ^d$, and $F$ is a tile of $\ZZ^d$. We also have $|\mathcal F_M|=3^{n^d/3}$.

To prove Theorem~\ref{thm:general-lower}, we show that for a random $M\in\mathcal M$ and a random $F\in\mathcal F_M$, we can asymptotically almost surely recover $M$ from $F$. That is to say, $F\notin \mathcal F_{M'}$ for any other $M'$. To show this, we will prove that for most $M$, the triangle matching $G_M$ is sufficiently different from all other $G_{M'}$. The claim below makes precise what sufficiently different means.

\begin{claim}\label{claim:compatible_prob}
Let $G,G'$ be perfect triangle matchings on the same $N$-vertex set, and let
$t\ge0$. If at least $t$ triangles of $G$ are not triangles of $G'$, then
they have at most $3^{N/3-ct}$ common maximal independent sets, where $c>0$
is an absolute constant.
\end{claim}
\begin{proof}
Without loss of generality, assume that there are exactly $t$ triangles of $G$ that are not triangles of $G'$.
Both graphs have $N/3$ triangles, so there are also exactly $t$ triangles of $G'$ that are not triangles of $G$. Consider the auxiliary graph $X$ whose vertices are the triangles of $G$
absent from $G'$, with an edge between two triangles when a triangle of $G'$ intersects both.
Observe that $\Delta(X)\le6$: each triangle of $G$ meets at most three
triangles of $G'$, each of which meets at most two other triangles of $G$.
Any graph with $t$ vertices and maximum degree $\Delta$ has an independent
set of size at least $t/(\Delta+1)$. Thus $X$ has an independent set
$\{T_1,\ldots,T_s\}$ with $s\ge |V(X)|/7=t/7$.

Choose a uniformly random maximal independent set $A$ of $G'$ by
independently selecting one vertex from each triangle. Since $T_i$ is
absent from $G'$, choose $a_i,b_i\in T_i$ lying in distinct triangles
of $G'$. Let $E_i$ be the event $\{a_i\in A,\ b_i\in A\}$, so $\PP(E_i)=1/9$.
No triangle of $G'$ intersects two of the $T_i$, so the events $E_i$
are independent. Each $E_i$ prevents $A$ from being independent in $G$.
Thus the number of common maximal independent sets is at most
\begin{align*}
    3^{N/3}\PP\left(\bigcap_{i=1}^s E_i^c\right)
    &=3^{N/3}(8/9)^s\\
    &\le3^{N/3}(8/9)^{t/7}
     =3^{N/3-ct}
\end{align*}
with $c=\frac{1}{7}\log_3 (9/8)$.
\end{proof}

Next, we show that for most $M$, $G_M$ is sufficiently different from all other $G_{M'}$. For $x\in\ZZ$, let $\operatorname{dist}(x,n\ZZ)=\min_{z\in\ZZ}|x-nz|$.

\begin{defn}
Let $t\ge1$ be an integer. A matrix $M\in\mathcal M$ is \emph{$t$-good} if
\[
    \operatorname{dist}\left(\sum_{i<j}r_{ij}M_{ij},n\ZZ\right)>t
\]
for every nonzero integer tuple $(r_{ij})_{i<j}$ satisfying $|r_{ij}|\le2^d$.
\end{defn}
\begin{claim}\label{claim:t-good-probability}
Fix $d\ge2$ and an integer $1\le t\le n/10$.
A uniformly random $M\in\mathcal M$ is $t$-good with probability at least
$1-O_d(t/n)$.
\end{claim}
\begin{proof}
Fix a nonzero integer tuple $(r_{ij})_{i<j}$ with
$|r_{ij}|\le2^d$, and put $g=\gcd(n,(r_{ij})_{i<j})\le2^d$.
The sum $\sum r_{ij}M_{ij}\pmod n$ is uniform on the multiples of $g$, so
\[
    \PP\!\left(\operatorname{dist}\left(\sum_{i<j}r_{ij}M_{ij},n\ZZ\right)
        \le t\right)\le\frac{2t/g+1}{n/g}
        =\frac{2t+g}{n}\le\frac{2t+2^d}{n}.
\]
There are $O_d(1)$ such tuples $(r_{ij})_{i<j}$, thus a union bound over the tuples proves the claim.
\end{proof}
\begin{claim}\label{claim:unique_triples}
Fix $d\ge2$ and an integer $1\le t\le n/10$.
If $M\in\mathcal M$ is $t$-good, then for every
$M'\in\mathcal M\setminus\{M\}$, at least $t$ triangles of $G_M$ are absent
from $G_{M'}$.
\end{claim}
\begin{proof}
Fix $M'\in\mathcal M\setminus\{M\}$. Let $i$ be the largest row index where
$M$ and $M'$ differ, and let $j$ be the smallest column index with
$M_{ij}\ne M'_{ij}$. First suppose $j<d$. With $e_j$ denoting the $j$-th
standard basis vector, put
\[
    X_0=\{se_j:0\le s<t\},\qquad
    X_1=\{(n-t+s)e_j:0\le s<t\},\qquad X=X_0\cup X_1.
\]
For each $x\in X$, let $y$ and $y'$ be the unique vertices with last
coordinate $n/3$ in the triangles containing $x$ in $G_M$ and $G_{M'}$,
respectively. We will show that $y\ne y'$ for every $x\in X_0$ or for
every $x\in X_1$. This implies that the $t$ triangles containing each $x\in X_0$ or $X_1$ in $G_M$ are absent from $G_{M'}$.

Write $v_k,v'_k$ for the $k$-th columns of $M,M'$. To construct $y$, add $v_d$
to $x$ and correct coordinates from $d-1$ upwards. Define integers
$c_{d-1},\ldots,c_1$, in this order, by setting $c_k$ to be the unique integer such that
\[
    0\le x_k+M_{kd}+\sum_{\ell=k+1}^{d-1}c_\ell M_{k\ell}+nc_k<n.
\]
Adding $c_kv_k$ corrects coordinate $k$ without disturbing the coordinates
below it, so $y=x+v_d+\sum_{k<d}c_kv_k$. The recursive definition gives
$|c_k|\le1+\sum_{\ell=k+1}^{d-1}|c_\ell|$, hence a simple induction gives
$|c_k|\le2^{d-k}\le2^d$. Define $c'_k$ analogously using $M'$.
Since rows below $i$ agree, for each $x\in X$ we have $c_k=c'_k$ for $k>i$.

As $x$ varies over $X$, only $x_j$ varies. The corrections below coordinate
$j$ do not depend on $x_j$, so all $c_k$ with $k>j$ remain fixed. Put
$B=M_{jd}+\sum_{k=j+1}^{d-1}c_kM_{jk}$. Since $|c_k|\le2^d$, $t$-goodness of $M$ gives
$\operatorname{dist}(B,n\ZZ)>t$. Let $c^*$ be the value of $c_j$ when $x_j=0$. Then $c^*$ is the unique integer such that $B+nc^*\in [0,n-1]$. Since $\operatorname{dist}(B,n\ZZ)>t$, we in fact have
\[
    t<B+nc^*<n-t.
\]
For $x\in X_0$, adding $0\le x_j<t$ keeps this coordinate in $[0,n-1]$,
so $c_j=c^*$. For $x\in X_1$, we have $n<x_j+B+nc^*<2n$, so we must
subtract one additional $n$, giving $c_j=c^*-1$.

Since $c_k=c'_k$ for $k>i$ and the entries in row $i$ before column $j$
agree, we obtain
\[
    y_i-y'_i\equiv c_j(M_{ij}-M'_{ij})+
    \sum_{k=j+1}^{d-1}c_k(M_{ik}-M'_{ik})+M_{id}-M'_{id}\pmod n.
\]
As $x$ varies over $X$, only $c_j$ varies on the right hand side and is constant on each of $X_0,X_1$. Let $\delta_l$ be the value of $y_i-y_i'\pmod{n}$ for any $x\in X_l$, $l=0,1$. Then $\delta_0-\delta_1\equiv M_{ij}-M'_{ij}\not\equiv0\pmod n$, thus one of $\delta_0,\delta_1$ is nonzero. Thus $y_i\ne y'_i$ throughout
at least one of $X_0,X_1$, so $y\ne y'$ as required.

If $j=d$, take $x=se_d$ for $0\le s<t$. The same correction formulas give
$y,y'\in Q$ with last coordinate $x_d+n/3$, since $t\le n/10$.
Again $c_k=c'_k$ for $k>i$, and the entries in row $i$ before column $d$
agree, so
$y_i-y'_i\equiv M_{id}-M'_{id}\not\equiv0\pmod n$.
Thus the triangles containing $x$ in $G_M$ and $G_{M'}$ differ for every
such $x$.
\end{proof}

\begin{proof}[Proof of Theorem~\ref{thm:general-lower}]
For $d=1$, there is just one matrix, and Observation~\ref{obs:lattice-triangles}
gives $3^{n/3}$ tiles. Suppose $d\ge2$, and put $t=\lfloor\sqrt n\rfloor$.
For sufficiently large $n$, we have $1\le t\le n/10$.
By Claim~\ref{claim:t-good-probability}, there are
$(1-o(1))n^{\binom d2}$ $t$-good matrices.
For each such $M$, Claims~\ref{claim:unique_triples}
and~\ref{claim:compatible_prob} give, by a union bound,
\[
    \left|\mathcal F_M\setminus
       \bigcup_{M'\ne M}\mathcal F_{M'}\right|
    \ge 3^{n^d/3}
       \left(1-n^{\binom d2}3^{-ct}\right) = (1-o(1))3^{n^d/3},
\]
with $c>0$ from Claim~\ref{claim:compatible_prob}. These families are pairwise disjoint. Summing over the $t$-good
matrices gives at least $(1-o(1))n^{\binom d2}3^{n^d/3}$ tiles, as required.
\end{proof}

\begin{remark}
The above construction depends on an ordering of the coordinates. By permuting the coordinates, one could obtain $T(n,d)\ge (1-o(1))d!n^{\binom d2}3^{n^d/3}$. However, this does not change the overall shape of $\Omega(n^{\binom d2}3^{n^d/3})$, so we omit the details for simplicity.
\end{remark}

\section{Concluding remarks}
\label{sec: conc}

Our results leave several questions concerning the enumeration and
structure of tiles. Throughout this section, random tiles are chosen
uniformly from the specified family.

Theorem~\ref{thm:one-dimensional} shows that almost all tiles in $[n]$
are AP tiles. A direction suggested to us by Nati Linial is to study
the exceptional family of non-AP tiles.

\begin{question}
What is the asymptotic growth of $T_{\mathrm{non\text{-}AP}}(n,1)$?
What is the typical structure of a random non-AP tile $A\subseteq[n]$,
and what is its typical density $|A|/n$?
\end{question}

When $3\mid n$, Theorems~\ref{thm:general-upper}
and~\ref{thm:general-lower} determine $T(n,d)$ up to a polynomial
factor in $n$ for every fixed $d$. Our lower bound with the factor
$n^{\binom d2}$ relies on the assumption $3\mid n$. Determining
the asymptotics in the other residue classes remains open in
dimensions $d\ge2$.

\begin{question}
For every fixed $d\ge2$, what is the asymptotic order of $T(n,d)$
as $n\to\infty$ through each of the residue classes $1$ and $2$
modulo $3$?
\end{question}

We call $A\subseteq\ZZ^d$ a
\emph{lattice tile} if $\ZZ^d=A\oplus L$ for some lattice
$L\subseteq\ZZ^d$. Our lower-bound construction in Theorem~\ref{thm:general-lower} consists of lattice
tiles, and Theorem~\ref{thm:one-dimensional} shows that almost all
one-dimensional tiles are lattice tiles.

\begin{question}
For every fixed $d\ge2$, is a random tile $A\subseteq[n]^d$ a lattice
tile with probability $1-o(1)$ as $n\to\infty$ through multiples
of $3$?
\end{question}

We end by highlighting that many of the open problems of \cite{BenjaminiKozmaTzalik2026} remain open. In particular, understanding tiles in domains other than $[n]^d$ remains an intriguing open problem.

\bibliographystyle{alpha}
\bibliography{ref}
\end{document}